\documentclass[11pt]{amsart}
\usepackage{mathrsfs}
\usepackage{amssymb,latexsym}
\usepackage{pstcol,pstricks,color}

 \numberwithin{equation}{section}

\newtheorem{theorem}{Theorem}[section]

\newtheorem{corollary}[theorem]{Corollary}

\newtheorem{lemma}[theorem]{Lemma}

\newtheorem{example}[theorem]{Example}

\def\qed{\hfill $\Box$}
\def\pf{\noindent {\it Proof.} }

\title{The largest singletons in weighted set partitions and its applications}

\begin{document}
\maketitle
\begin{center}
Yidong Sun$^\dag$ and Yanjie Xu


Department of Mathematics, Dalian Maritime University, 116026 Dalian, P.R. China\\[5pt]

{\it $^\dag$Email: sydmath@yahoo.com.cn,  }
\end{center}\vskip0.5cm

\subsection*{Abstract}
Recently, Deutsch and Elizalde studied the largest and the smallest fixed points of
permutations. Motivated by their work, we consider the analogous problems in weighted set partitions.
Let $A_{n,k}(\mathbf{t})$ denote the total weight of partitions on $[n+1]$ with the largest singleton $\{k+1\}$.
In this paper, explicit formulas for $A_{n,k}(\mathbf{t})$ and many combinatorial identities involving $A_{n,k}(\mathbf{t})$
are obtained by umbral operators and combinatorial methods. As applications, we investigate three special cases
such as permutations, involutions and labeled forests. Particularly in the permutation case, we derive a surprising identity analogous
to the Riordan identity related to tree enumerations, namely,
\begin{eqnarray*}
\sum_{k=0}^{n}\binom{n}{k}D_{k+1}(n+1)^{n-k}  &=& n^{n+1},
\end{eqnarray*}
where $D_{k}$ is the $k$-th derangement number or the number of permutations of $\{1,2,\dots, k\}$ with no fixed points.

\medskip

{\bf Keywords}: Set partition; Bell polynomial; Permutation; Involution; Labeled forest.

\noindent {\sc 2000 Mathematics Subject Classification}: Primary
05A05; Secondary 05C30

\section{Introduction}

A {\it partition} of a set $[n]=\{1, 2, \dots, n\}$ is a
collection $\pi=\{\mathbb{B}_1, \mathbb{B}_2, \dots, \mathbb{B}_r\}$ of nonempty and mutually disjoint subsets of $[n]$,
called {\it blocks}, whose union is $[n]$. For a block $\mathbb{B}$, we denote by $|\mathbb{B}|$ the size of the block $\mathbb{B}$, that is the
number of the elements in the block $\mathbb{B}$. A block $\mathbb{B}$ will be called {\it singleton} if $|\mathbb{B}| = 1$.
If $\{k\}$ is a singleton of a partition, we denote it by $k$ for short. If $|\mathbb{B}| = j$, we assign a weight $t_j$ for $\mathbb{B}$. The
weight $w(\pi)$ of a partition $\pi$ is defined to be the product of the weight of each block of $\pi$.

It is well known that the weight of partitions of $[n]$ with $r$ blocks is the
partial Bell polynomial ${\mathcal{B}}_{n,r}\big(t_1,t_2,\dots\big)$ \cite{Comtet} on the variables $\{t_j\}_{j\geq 1}$, that is
\begin{eqnarray*}
{\mathcal{B}}_{n,r}\big(t_1,t_2,\dots\big)&=&\sum_{\kappa_n(r)}\frac{n!}{r_1!r_2!\cdots
r_n!}\left(\frac{t_1}{1!}\right)^{r_1}\left(\frac{t_2}{2!}\right)^{r_2}\cdots
\left(\frac{t_n}{n!}\right)^{r_n},
\end{eqnarray*}
where the summation $\kappa_n(r)$ is for all the nonnegative integer solutions of $r_1+r_2+\cdots+r_n=r$ and $r_1+2r_2+\cdots +nr_n=n$.
And the total weight for partitions of $[n]$ is the complete Bell polynomial
\begin{eqnarray*}
{\mathcal{Y}}_n(\mathbf{t})={\mathcal{Y}}_{n}\big(t_1,t_2,\dots\big)=\sum_{r=0}^{n}{\mathcal{B}}_{n,r}\big(t_1,t_2,\dots\big),
\end{eqnarray*}
which has the exponential generating function
\begin{eqnarray*}
{\mathcal{Y}}(\mathbf{t};x)=\sum_{n\geq 0}{\mathcal{Y}}_{n}\big(t_1,t_2,\dots\big)\frac{x^n}{n!}=\exp{\big(\sum_{j\geq 1}t_j\frac{x^j}{j!}\big)}.
\end{eqnarray*}

Let $\mathbb{A}_{n,k}$ denote the set of partitions of $[n+1]$ with
the largest singleton $k+1$. Let $A_{n,k}(\mathbf{t})$ denote
the total weight of partitions in $\mathbb{A}_{n,k}$. Clearly,
\begin{eqnarray*}
A_{n,0}(\mathbf{t})=t_1{\mathcal{Y}}_{n}\big(0,t_2,\dots\big) \hskip.5cm {\mbox{and}} \hskip.5cm A_{n,n}(\mathbf{t})=t_1{\mathcal{Y}}_{n}\big(t_1,t_2,\dots\big),
\end{eqnarray*}
where ${\mathcal{Y}}_{n}\big(0,t_2,\dots\big)$ is the weight of partitions of $[n]$ without singletons.

Recently, Deutsch and Elizalde \cite{DeuEliz} studied the largest fixed points of
permutations, which is the special case when $t_j=(j-1)!$ for $j\geq 1$. Later, Sun and Wu \cite{SunWu} considered the
largest singletons in set partitions, which is the special case when $t_j=1$ for $j\geq 1$.

In this paper we will investigate the largest singletons in weighted set partitions generally. The next section is devoted to
studying the properties of $A_{n,k}(\mathbf{t})$, involving its explicit formulas and many combinatorial identities for
$A_{n,k}(\mathbf{t})$. In the third section, we consider the permutation case, i.e., the special case when $t_j=(j-1)!$ for $j\geq 1$,
and derive a surprising identity analogous to the Riordan identity related to tree enumerations. In the forth section,
we study the involution case which is the special case when $t_1=t_2=1, t_j=0$ for $j\geq 3$. In the final section,
we focus on the labeled forest case which is the special case when $t_j=j^{j-1}$ for $j\geq 1$.

\section{The properties of $A_{n,k}(\mathbf{t})$}

According to the definition of $A_{n,k}(\mathbf{t})$, for any weighted partition $\pi$ of
$[n+1]$ with the largest singleton $k+1$, if $k$ is also a singleton, delete the singleton $k+1$ and subtracting one from all the entries
large than $k+1$, we obtain a partition of $[n]$ with the largest singleton $k$.
This contributes the weight $t_1A_{n-1,k-1}(\mathbf{t})$; if $k$ is not a singleton, exchange $k$ and $k+1$,
we obtain a partition of $[n+1]$ with the largest singleton $k$. This contributes the weight $A_{n,k-1}(\mathbf{t})$. Then
we obtain a recurrence for $n, k\geq 1$,
\begin{eqnarray}\label{eqn 2.1}
A_{n,k}(\mathbf{t})=A_{n,k-1}(\mathbf{t})+t_1A_{n-1, k-1}(\mathbf{t})
\end{eqnarray}
with the initial conditions $A_{n,0}(\mathbf{t})=t_1{\mathcal{Y}}_{n}\big(0,t_2,\dots\big)$ for $n\geq 0$.

\begin{lemma}\label{lemma 2.1}
The bivariate exponential generating function for $A_{n+k,k}(\mathbf{t})$ is
given by
\begin{eqnarray*}
A(\mathbf{t}; x, y)=\sum_{n,k\geq
0}A_{n+k,k}(\mathbf{t})\frac{x^n}{n!}\frac{y^k}{k!}=t_1e^{-xt_1}{\mathcal{Y}}(\mathbf{t};x+y).
\end{eqnarray*}
\end{lemma}

\pf Define
\begin{eqnarray*}
A_k(\mathbf{t}; x)=\sum_{n\geq 0}A_{n+k,k}(\mathbf{t})\frac{x^n}{n!}.
\end{eqnarray*}
Clearly, $A_0(\mathbf{t}; x)=t_1e^{-xt_1}{\mathcal{Y}}(\mathbf{t}; x)$. From
(\ref{eqn 2.1}), one can derive that
\begin{eqnarray*}
A_k(\mathbf{t}; x)=t_1A_{k-1}(\mathbf{t}; x)+\frac{\partial}{\partial x}A_{k-1}(\mathbf{t}; x),
\end{eqnarray*}
which produces
\begin{eqnarray*}
A_k(\mathbf{t}; x)=(t_1+\frac{\partial}{\partial x})A_{k-1}(\mathbf{t}; x)=(t_1+\frac{\partial}{\partial x})^kA_0(\mathbf{t}; x).
\end{eqnarray*}
Then
\begin{eqnarray*}
A(\mathbf{t}; x,y) &=& \sum_{k\geq 0}A_k(\mathbf{t}; x)\frac{y^k}{k!}=\sum_{k\geq 0}\frac{y^k(t_1+\frac{\partial}{\partial x})^k}{k!}A_0(\mathbf{t}; x) \\
                   &=& e^{yt_1+y\frac{\partial}{\partial x}}t_1e^{-xt_1}{\mathcal{Y}}(\mathbf{t}; x)=t_1e^{yt_1}e^{y\frac{\partial}{\partial x}}e^{-xt_1}{\mathcal{Y}}(\mathbf{t}; x) \\
                   &=& t_1e^{yt_1}e^{-(x+y)t_1}{\mathcal{Y}}(\mathbf{t}; x+y)=t_1e^{-xt_1}{\mathcal{Y}}(\mathbf{t}; x+y).
\end{eqnarray*}
This completes the proof. \qed\vskip.2cm

\begin{theorem} For any integers $n,m\geq 0$ and any indeterminant $\lambda$, there hold
\begin{eqnarray}
\sum_{k=0}^{n}\binom{k+\lambda-1}{k}A_{n+m,m+k}(\mathbf{t}) &=& \sum_{k=0}^{n}\binom{n+\lambda}{k}\binom{n+\lambda-k-1}{n-k}A_{m+k,m}(\mathbf{t})t_1^{n-k},\label{eqn 2.2A}\\
\sum_{k=0}^{n}\binom{k+\lambda-1}{k}A_{n+m,m+k}(\mathbf{t}) &=& \sum_{k=0}^{n}(-1)^{n-k}\binom{n+\lambda}{k}\mathcal{Y}_{m+k}(\mathbf{t})t_1^{n-k+1}. \label{eqn 2.2B}
\end{eqnarray}
\end{theorem}
\pf With the umbra $\mathbf{Y_t}$, given by $\mathbf{Y_t}^{n}=\mathcal{Y}_{n}(\mathbf{t})$, $\mathcal{Y}(\mathbf{t}; x)$ may be written as
${\mathcal{Y}}(\mathbf{t}; x)=e^{\mathbf{Y_t}x}$. (See \cite{Gesselb, Roman, RomRota} for more information on umbral calculus, to cite only a few).
Then, by Lemma 2.1, we have
\begin{eqnarray*}
A(\mathbf{t}; x,y)=t_1e^{\mathbf{Y_t}(x+y)-t_1x}=t_1e^{(\mathbf{Y_t}-t_1)x}e^{\mathbf{Y_t}y}.
\end{eqnarray*}
When comparing the coefficient of $\frac{x^ny^k}{n!k!}$, $A_{n+k,k}(\mathbf{t})$ can be represented umbrally as
\begin{eqnarray}\label{eqn 2.2C}
A_{n+k,k}(\mathbf{t})=t_1\mathbf{Y_t}^{k}(\mathbf{Y_t}-t_1)^{n}.
\end{eqnarray}
Let $[x^n]f(x)$ denote the coefficient of $x^n$ in the formal power series $f(x)$, we get
\begin{eqnarray*}
\lefteqn{\sum_{k=0}^{n}\binom{k+\lambda-1}{k}A_{n+m,m+k}(\mathbf{t})}    \\
&=& \sum_{k=0}^{n}(-1)^{k}\binom{-\lambda}{k}t_1\mathbf{Y_t}^{m+k}(\mathbf{Y_t}-t_1)^{n-k}         \\
&=& t_1\mathbf{Y_t}^{m}(\mathbf{Y_t}-t_1)^{n}\sum_{k=0}^{n}\binom{-\lambda}{k}\Big(-\frac{\mathbf{Y_t}}{\mathbf{Y_t}-t_1}\Big)^{k}    \\
&=& t_1\mathbf{Y_t}^{m}(\mathbf{Y_t}-t_1)^{n}\sum_{k=0}^{n}[x^k]\Big(1-\frac{x\mathbf{Y_t}}{\mathbf{Y_t}-t_1}\Big)^{-\lambda}    \\
&=& t_1\mathbf{Y_t}^{m}(\mathbf{Y_t}-t_1)^{n}[x^n]\frac{1}{1-x}\Big(1-\frac{x\mathbf{Y_t}}{\mathbf{Y_t}-t_1}\Big)^{-\lambda}    \\
&=& t_1\mathbf{Y_t}^{m}(\mathbf{Y_t}-t_1)^{n}[x^n]\frac{1}{(1-x)^{\lambda+1}}\Big(1-\frac{x}{(1-x)}\frac{t_1}{(\mathbf{Y_t}-t_1)}\Big)^{-\lambda}    \\
&=& t_1\mathbf{Y_t}^{m}(\mathbf{Y_t}-t_1)^{n}[x^n]\sum_{k=0}^{n}\binom{-\lambda}{n-k}\frac{x^{n-k}}{(1-x)^{n+\lambda-k+1}}
\Big(-\frac{t_1}{\mathbf{Y_t}-t_1}\Big)^{n-k}    \\
&=& \sum_{k=0}^{n}(-1)^k\binom{-(n+\lambda-k+1)}{k}\binom{-\lambda}{n-k}t_1\mathbf{Y_t}^{m}(\mathbf{Y_t}-t_1)^{k}(-t_1)^{n-k}   \\
&=& \sum_{k=0}^{n}\binom{n+\lambda}{k}\binom{n+\lambda-k-1}{n-k}A_{m+k,m}(\mathbf{t})t_1^{n-k},   \\
\end{eqnarray*}
which proves (\ref{eqn 2.2A}).

By the identity
\begin{eqnarray*}
\binom{n}{k}\binom{k}{i}=\binom{n}{i}\binom{n-i}{k-i},
\end{eqnarray*}
and the Vandermonde's convolution identity
\begin{eqnarray*}
\sum_{k=0}^{n}\binom{a}{k}\binom{b}{n-k}=\binom{a+b}{n},
\end{eqnarray*}
we have
\begin{eqnarray*}
\lefteqn{\sum_{k=0}^{n}\binom{k+\lambda-1}{k}A_{n+m,m+k}(\mathbf{t})}    \\
&=& \sum_{k=0}^{n}\binom{n+\lambda}{k}\binom{-\lambda}{n-k}t_1\mathbf{Y_t}^{m}(\mathbf{Y_t}-t_1)^{k}(-t_1)^{n-k}   \\
&=& \sum_{k=0}^{n}\binom{n+\lambda}{k}\binom{-\lambda}{n-k}\sum_{i=0}^{k}\binom{k}{i}t_1\mathbf{Y_t}^{m+i}(-t_1)^{n-i}  \\
&=& \sum_{i=0}^{n}t_1\mathbf{Y_t}^{m+i}(-t_1)^{n-i}\sum_{k=i}^{n}\binom{n+\lambda}{k}\binom{-\lambda}{n-k}\binom{k}{i}  \\
&=& \sum_{i=0}^{n}\binom{n+\lambda}{i}t_1\mathbf{Y_t}^{m+i}(-t_1)^{n-i}\sum_{k=i}^{n}\binom{-\lambda}{n-k}\binom{n+\lambda-i}{k-i}  \\
&=& \sum_{i=0}^{n}\binom{n+\lambda}{i}t_1\mathbf{Y_t}^{m+i}(-t_1)^{n-i} \\
&=& \sum_{k=0}^{n}(-1)^{n-k}\binom{n+\lambda}{k}\mathcal{Y}_{m+k}(\mathbf{t})t_1^{n-k+1},
\end{eqnarray*}
which proves (\ref{eqn 2.2B}). \qed \vskip0.2cm

The case $\lambda=0$ in (\ref{eqn 2.2B}), yields the explicit formula for $A_{n+m,m}(\mathbf{t})$.
\begin{corollary} For any integers $n, m\geq 0$, there holds
\begin{eqnarray}\label{eqn 2.3A}
A_{n+m,m}(\mathbf{t})=\sum_{k=0}^{n}(-1)^{n-k}\binom{n}{k}t_1^{n-k+1}{\mathcal{Y}}_{m+k}(\mathbf{t}).
\end{eqnarray}
\end{corollary}
\pf Let $\mathbb{X}$ denote the set of partitions of $[n+m+1]$ containing at least the singleton $m+1$.
Clearly, $\mathbb{X}$ has the weight $t_1\mathcal{Y}_{n+m}(\mathbf{t})$. Let $\mathbb{X}_i$ be
the subset of $\mathbb{X}$ containing another singleton $m+i+1$ for $1\leq i\leq n$. Set $\overline{\mathbb{X}}_i=\mathbb{X}-\mathbb{X}_i$,
then $\bigcap_{i=1}^n\overline{\mathbb{X}}_i$ is just the set of partitions of $[n+m+1]$ with the largest singleton $m+1$,
so $\bigcap_{i=1}^n\overline{\mathbb{X}}_i$ has the weight $A_{n+m,m}(\mathbf{t})$.
For any nonempty $(n-k)$-subset $\mathbb{S}\in [n]$, $\bigcap_{i\in \mathbb{S}}\mathbb{X}_i$
is the set of partitions of $[n+m+1]$ containing at least the number $n-k+1$ of singletons $m+1$ and $m+i+1$ for all $i\in \mathbb{S}$,
so $\bigcap_{i\in \mathbb{S}}\mathbb{X}_i$ has the weight $t_1^{n-k+1}{\mathcal{Y}}_{m+k}(\mathbf{t})$.
By the Inclusion-Exclusion principle, we have
\begin{eqnarray*}
w(\bigcap_{i=1}^n\overline{\mathbb{X}}_i) &=& w(\mathbb{X}-\bigcup_{i=1}^{n}\mathbb{X}_i) \\
                                          &=& w(\mathbb{X})+\sum_{k=0}^{n-1}(-1)^{n-k}\binom{n}{k}w(\bigcap_{i\in \mathbb{S}, |\mathbb{S}|=n-k}\mathbb{X}_i) \\
                                          &=& \sum_{k=0}^{n}(-1)^{n-k}\binom{n}{k}t_1^{n-k+1}{\mathcal{Y}}_{m+k}(\mathbf{t}),
\end{eqnarray*}
which proves (\ref{eqn 2.3A}). \qed \vskip0.2cm

\begin{corollary} For any integers $n,m\geq 0$, there holds
\begin{eqnarray}
\sum_{k=0}^{n}A_{n+m,m+k}(\mathbf{t}) &=& \frac{t_1\mathcal{Y}_{n+m+1}(\mathbf{t})-A_{n+m+1,m}(\mathbf{t})}{t_1}, \label{eqn 2.3B}   \\
\sum_{k=0}^{n}(k+1)A_{n+m,m+k}(\mathbf{t}) &=& \frac{A_{n+m+2,m}(\mathbf{t})-t_1\mathcal{Y}_{n+m+2}(\mathbf{t})
+(n+2)t_1^{2}\mathcal{Y}_{n+m+1}(\mathbf{t})}{t_1^2},   \label{eqn 2.3C}    \\
\sum_{k=0}^{n}(n-k+1)A_{n+m,m+k}(\mathbf{t}) &=& \frac{t_1\mathcal{Y}_{n+m+2}(\mathbf{t})-A_{n+m+2,m}(\mathbf{t})-(n+2)t_1A_{n+m+1,m}(\mathbf{t})
}{t_1^2}.  \label{eqn 2.3D}
\end{eqnarray}
\end{corollary}
\pf The case $n:=n+1$ in (\ref{eqn 2.3A}), together with the case $\lambda=1$ in (\ref{eqn 2.2B}), yields (\ref{eqn 2.3B}).
The case $n:=n+2$ in (\ref{eqn 2.3A}), together with the case $\lambda=2$ in (\ref{eqn 2.2B}), yields (\ref{eqn 2.3C}). And (\ref{eqn 2.3D}) can be
easily obtained from (\ref{eqn 2.3B}) and (\ref{eqn 2.3C}). \qed \vskip0.2cm

\begin{theorem} For any integers $n,m,k\geq 0$, there holds
\begin{eqnarray}
A_{n+m+k,m+k}(\mathbf{t})&=& \sum_{j=0}^{m}\binom{m}{j}t_1^{m-j}A_{n+k+j,k}(\mathbf{t}). \label{eqn 2.4}
\end{eqnarray}
\end{theorem}
\pf Here we provide a combinatorial proof. For any $\pi\in \mathbb{A}_{n+m+k,m+k}$,
suppose that $\pi$ has exactly $m-j$ singletons in $\{k+1, \dots, k+m\}$ which contribute the weight $t_1^{m-j}$,
and there are $\binom{m}{j}$ ways to do this. The remainder $j$
elements in $\{k+1, \dots, k+m\}$ can not be singletons in $\pi$. These $j$ elements can be regarded as the roles that greater than $m+k+1$,
so the remainder $n+k+j+1$ elements can be partitioned with the largest singleton $m+k+1$, which contributes the weight $A_{n+k+j,k}(\mathbf{t})$.
Thus the total weight of such partitions is $\binom{m}{j}t_1^{m-j}A_{n+k+j,k}(\mathbf{t})$.
Summing up all the possible cases yields (\ref{eqn 2.4}).  \qed\vskip.2cm

\begin{theorem} For any integers $n,m\geq 0$ and any indeterminant $y$, there hold
\begin{eqnarray}
\sum_{k=0}^{n}\binom{n}{k}A_{n+m,m+k}(\mathbf{t})y^{k}  &=& \sum_{k=0}^{n}(-1)^{n-k}\binom{n}{k}\mathcal{Y}_{m+k}(\mathbf{t})(y+1)^{k}t_1^{n-k+1}, \label{eqn 2.5} \\
\sum_{k=0}^{n}\binom{n}{k}A_{m+k,m}(\mathbf{t})y^{n-k}  &=& t_1\sum_{k=0}^{n}\binom{n}{k}\mathcal{Y}_{m+k}(\mathbf{t})(y-t_1)^{n-k}.  \label{eqn 2.6}
\end{eqnarray}
\end{theorem}
\pf By (\ref{eqn 2.2C}), we have
\begin{eqnarray*}
\sum_{k=0}^{n}\binom{n}{k}A_{n+m,m+k}(\mathbf{t})y^{k} &=& \sum_{k=0}^{n}\binom{n}{k}t_1\mathbf{Y_t}^{m+k}(\mathbf{Y_t}-t_1)^{n-k}y^{k}  \\
                                                       &=& t_1\mathbf{Y_t}^{m}((y+1)\mathbf{Y_t}-t_1)^{n}  \\
                                                       &=& \sum_{k=0}^{n}(-1)^{n-k}\binom{n}{k}(y+1)^{k}\mathbf{Y_t}^{m+k}t_1^{n-k+1} \\
                                                       &=& \sum_{k=0}^{n}(-1)^{n-k}\binom{n}{k}(y+1)^{k}\mathcal{Y}_{m+k}(\mathbf{t})t_1^{n-k+1}, \\
\end{eqnarray*}
which proves (\ref{eqn 2.5}). Similarly, (\ref{eqn 2.6}) can be obtained, but here we provide a combinatorial proof.

Let $\mathbb{X}_{n,m}=\bigcup_{k=0}^n\mathbb{X}_{n, m, k}$ and $\mathbb{X}_{n, m, k}$ denote the set of pairs $(\pi, \mathbb{S})$ such that
\begin{itemize}
\item $\mathbb{S}$ is an $(n-k)$-subset of $[m+2, n+m+1]=\{m+2, \dots, n+m+1\}$, and each element of $\mathbb{S}$ is colored by $t_1$ or $y-t_1$;
\item $\pi$ is a partition of the set $[n+m+1]-\mathbb{S}$ with the largest singleton $m+1$, and each element of $[n+m+1]-\mathbb{S}$ is only colored by $1$.
\end{itemize}

Let $\mathbb{Y}_{n,m}=\bigcup_{k=0}^n\mathbb{Y}_{n, m, k}$ and $\mathbb{Y}_{n, m, k}$ denote the set of pairs $(\pi, \mathbb{S})$ such that
\begin{itemize}
\item $\mathbb{S}$ is an $(n-k)$-subset of $[m+2, n+m+1]$ and each element of $\mathbb{S}$ is only colored by $y-t_1$;
\item $\pi$ is a partition of the set $[n+m+1]-\mathbb{S}$ such that $m+1$ must be a singleton, and each element of $[n+m+1]-\mathbb{S}$ is only colored by $1$.
\end{itemize}
The weight of $(\pi, \mathbb{S})$ is defined to be the product of the weight of $\pi$ and the color of each element of $[n+m+1]$. Clearly, the weights of $\mathbb{X}_{n,m}$
and $\mathbb{Y}_{n,m}$ are counted respectively by the left and right sides of (\ref{eqn 2.6}).

Given any pair $(\pi, \mathbb{S})\in \mathbb{X}_{n,m}$, $\mathbb{S}$ can be partitioned into two parts $\mathbb{S}_1$ and $\mathbb{S}_{2}$ such that each
element of $\mathbb{S}_1$ is colored by $y-t_1$ and each element of $\mathbb{S}_2$ is colored by $t_1$.
Regard each element of $\mathbb{S}_2$ as a singleton which is weighted by $t_1$ and colored by $1$,
together with $\pi$, we obtain a partition $\pi_1$ of $[n+m+1]-\mathbb{S}_1$ such that $m+1$ is always a singleton.
Then the pair $(\pi_1, \mathbb{S}_1)$ lies in $\mathbb{Y}_{n,m}$.

Conversely, for any pair $(\pi_1, \mathbb{S}_1)\in \mathbb{Y}_{n,m}$, let $\mathbb{S}$ denote the union of $\mathbb{S}_1$ and the singletons of $\pi_1$
greater than $m+1$, then $\pi_1$ can be partitioned into two parts $\pi$ and $\pi'$ such that $\pi$ is a partition of $[n+m+1]-\mathbb{S}$ with
the largest singleton $m+1$ and $\pi'$ is the singletons of $\pi_1$ greater than $m+1$. Regard $\pi'$ as a subset of $[m+2, n+m+1]$ in which each element
is colored by $t_1$, together with $\mathbb{S}_1$, we obtain an $(n-k)$-subset of $[m+2, n+m+1]$ for some $k$ such that
each element of $\mathbb{S}$ is colored by $t_1$ or $y-t_1$. Then the pair $(\pi, \mathbb{S})$ lies in $\mathbb{X}_{n,m}$.

Clearly we find a bijection between $\mathbb{X}_{n,m}$ and $\mathbb{Y}_{n,m}$, which proves (\ref{eqn 2.6}). \qed\vskip.2cm

The cases $y=-1$ in (\ref{eqn 2.5}) and $y=t_1$ in (\ref{eqn 2.6}) lead to
\begin{corollary} For any integers $n, m \geq 0$, there hold
\begin{eqnarray*}
\sum_{k=0}^{n}(-1)^{n-k}\binom{n}{k}A_{n+m,m+k}(\mathbf{t}) &=& \mathcal{Y}_{m}(\mathbf{t})t_1^{n+1},  \\
\sum_{k=0}^{n}\binom{n}{k}A_{m+k,m}(\mathbf{t})t_1^{n-k-1}  &=& \mathcal{Y}_{m+n}(\mathbf{t}).
\end{eqnarray*}
\end{corollary}

The case $y:=\frac{yt_1}{y+1}$ in (\ref{eqn 2.6}), together with (\ref{eqn 2.5}) generates the following result which has a combinatorial interpretation.
\begin{corollary} For any integers $n, m \geq 0$, there holds
\begin{eqnarray}\label{eqn 2.7}
\sum_{k=0}^{n}\binom{n}{k}A_{m+k,m}(\mathbf{t})(y+1)^{k}(yt_1)^{n-k} &=& \sum_{k=0}^{n}\binom{n}{k}A_{n+m,m+k}(\mathbf{t})y^{k}.
\end{eqnarray}
\end{corollary}
\pf Let $\mathbb{X}_{n,m}^{*}=\bigcup_{k=0}^n\mathbb{X}_{n, m, k}^{*}$ and $\mathbb{X}_{n, m, k}^{*}$ denote the set of pairs $(\pi, \mathbb{S})$ such that
\begin{itemize}
\item $\pi$ is a partition of the set $[n+m+1]$ containing at least the singleton $m+1$;
\item $\mathbb{S}$ is an $(n-k)$-subset of $[m+2, n+m+1]$ which is also the set of singletons of $\pi$ greater than $m+1$,
each element of $\mathbb{S}$ is only colored by $y$ and each element of $[m+2, n+m+1]-\mathbb{S}$ is colored by $1$ or $y$;
\item each element of $[m+1]$ is only colored by $1$.

\end{itemize}

Let $\mathbb{Y}_{n,m}^{*}=\bigcup_{k=0}^n\mathbb{Y}_{n, m, k}^{*}$ and $\mathbb{Y}_{n, m, k}^{*}$ denote the set of pairs $(\pi, \mathbb{S})$ such that
\begin{itemize}
\item $\mathbb{S}$ is a $k$-subset $\{i_1, i_2, \dots, i_k\}$ of $[m+2, n+m+1]$ in increasing order, each element of $\mathbb{S}$ is only colored
by $y$ and each element of $[n+m+1]-\mathbb{S}$ is only colored by $1$;
\item $\pi$ is a partition of the set $[n+m+1]$ such that $i_k$ must be the largest singleton if $\mathbb{S}$ is not empty and
$m+1$ must be the largest singleton if $\mathbb{S}$ is empty;
\item each element of $[m+2,n+m+1]-\mathbb{S}$ must not be a singleton.
\end{itemize}
The weight of $(\pi, \mathbb{S})$ is defined to be the product of the weight of $\pi$ and the colors of all elements in $[n+m+1]$. Clearly, any
$(\pi, \mathbb{S})\in \mathbb{X}_{n,m}^{*}$ can be obtained as follows. First choose an $(n-k)$-subset $\mathbb{S}$ of $[m+2, n+m+1]$, there are
$\binom{n}{k}$ ways to do this. Regard each element of $\mathbb{S}$ as a singleton with color $y$. Then color each element of $[m+2, n+m+1]-\mathbb{S}$
by $1$ or $y$, namely, each element of $[m+2, n+m+1]-\mathbb{S}$ is colored by $y+1$. Now partitioning $[n+m+1]-\mathbb{S}$ such that the largest singleton is
$m+1$, together with the $n-k$ singletons formed form $\mathbb{S}$, we get the partition $\pi$ of $[n+m+1]$ such that $m+1$ must be a singleton;
Hence the total weight of pairs $(\pi, \mathbb{S})\in \mathbb{X}_{n,m}^{*}$ is just the left hand side of (\ref{eqn 2.7}).

Similarly, the total weight of pairs $(\pi, \mathbb{S})\in \mathbb{Y}_{n,m}^{*}$ is just the right hand side of (\ref{eqn 2.7}) if
regarding each element of $[m+2,n+m+1]-\mathbb{S}$ as the role greater than $i_k$ when $\mathbb{S}$ is not empty.

Now we can construct a bijection $\varphi$ between $\mathbb{X}_{n,m}^{*}$ and $\mathbb{Y}_{n,m}^{*}$ which preserves the weights.
For any $(\pi, \mathbb{S})\in \mathbb{X}_{n,m}^{*}$, let $\mathbb{S}_1$ denote the set of elements of $[n+m+1]$ with colors $y$. Clearly, $\mathbb{S}$
is a subset of $\mathbb{S}_1$. Assume that $\mathbb{S}_1=\{i_1, i_2, \dots, i_k\}$ for some $0\leq k\leq n$ in increasing order. If $\mathbb{S}_1$ is the
empty set $\emptyset$, which implies that $\mathbb{S}=\emptyset$ and all elements of $[n+m+1]$ are colored by $1$,
it is obvious that $(\pi, \emptyset)\in \mathbb{Y}_{n,m}^{*}$. Then
define $\varphi(\pi, \emptyset)=(\pi, \emptyset)$. If $\mathbb{S}_1$ is not the empty set, exchanging $m+1$ and $i_k$ in $\pi$,
we obtain a partition $\pi_1$, it is easily to verify that $(\pi_1, \mathbb{S}_1)\in \mathbb{Y}_{n,m}^{*}$ and has the same weight
as $(\pi, \mathbb{S})$. Then define $\varphi(\pi, \mathbb{S})=(\pi_1, \mathbb{S}_1)$.

Conversely, for any $(\pi_1, \mathbb{S}_1)\in \mathbb{Y}_{n,m}^{*}$, if $\mathbb{S}_1=\emptyset$, so $\pi_1$ has the largest singleton $m+1$,
then $(\pi_1, \emptyset)\in \mathbb{X}_{n,m}^{*}$ and define $\varphi^{-1}(\pi_1, \emptyset)=(\pi_1, \emptyset)$. If $\mathbb{S}_1\neq\emptyset$,
assume that $\mathbb{S}_1=\{i_1, i_2, \dots, i_k\}$ for some $1\leq k\leq n$ in increasing order,
let $\mathbb{S}$ denote the set of all the elements in $\mathbb{S}_1$ such that each forms a singleton of $\pi_1$. Now exchanging $m+1$ and $i_k$
in $\pi_1$, we obtain a partition $\pi$, it is easy verifiable that $(\pi, \mathbb{S})\in \mathbb{X}_{n,m}^{*}$ which has the same weight as
$(\pi_1, \mathbb{S}_1)$. Then define $\varphi^{-1}(\pi_1, \mathbb{S}_1)=(\pi, \mathbb{S})$.

Clearly, $\varphi$ is indeed a bijection between $\mathbb{X}_{n,m}^{*}$ and $\mathbb{Y}_{n,m}^{*}$, which proves (\ref{eqn 2.7}). \qed\vskip.2cm

\section{The special case for permutations }

In this section, we consider the special case when $t_j=(j-1)!$ for $j\geq 1$. That is to assign a cycle structure to each block of partitions of $[n+1]$,
such partitions with weight $\mathbf{t}=(0!, 1!, 2!, \dots)$ is equivalent to permutations of $[n+1]$. Let $P_{n,k}=A_{n,k}(\mathbf{t})$ with
$\mathbf{t}=(0!, 1!, 2!, \dots)$, namely, $P_{n,k}$ is the number of permutations of $[n+1]$ with the largest fixed point $k+1$.
From (\ref{eqn 2.3A}) and (\ref{eqn 2.4}), one has the explicit formulas for $P_{n,k}$
\begin{eqnarray*}
P_{n+k,k}=\sum_{j=0}^{n}(-1)^{n-j}\binom{n}{j}(k+j)!=\sum_{j=0}^{k}\binom{k}{j}D_{n+j}.
\end{eqnarray*}
Clearly, $P_{n,n}=n!=\mathcal{Y}_n(0!, 1!, 2!, \dots)$ and $P_{n,0}=D_n=\mathcal{Y}_n(0, 1!, 2!, \dots)$,
where $D_n$ is the derangement number of $[n]$, i.e., the number of permutations of $[n]$ without fixed points.
See Table $1$ for some small values of $P_{n,k}$.

\begin{center}
\begin{eqnarray*}
\begin{array}{c|ccccccccc}\hline
n/k  &   0   &   1   &   2   &  3    &  4    &  5     &  6       \\\hline
  0  &   1   &       &       &       &       &        &          \\
  1  &   0   &   1   &       &       &       &        &           \\
  2  &   1   &   1   &   2   &       &       &        &            \\
  3  &   2   &   3   &   4   &   6   &       &        &           \\
  4  &   9   &  11   &  14   &  18   &  24   &        &           \\
  5  &  44   &  53   &  64   &  78   &  96   &  120   &              \\
  6  & 265   & 309   & 362   & 426   & 504   &  600   & 720         \\\hline
\end{array}
\end{eqnarray*}
Table 1. The values of $P_{n,k}$ for $n$ and $k$ up to $6$.
\end{center}

In fact $\{P_{n,k}\}_{n\geq k\geq 0}$ forms the difference table introduced by Euler, which has been
investigated in depth in the derangement theory \cite{ClarkHanZeng, Dumont, DumontRand, RakotonA, RakotonB}. Chen \cite{Chen} also gave another
two interpretations for $P_{n,k}$ using $k$-relative derangements on $[n]$ and skew derangements from $[n]$
to $\{-k+1, \dots, -1, 0, 1, \dots, n-k\}$ for $0\leq k\leq n$. Actually, Chen established a bijection between these two settings.
In a forthcoming paper, we find the bijective connections between several combinatorial objects which are counted by the Euler difference table.
Recently, Deutsch and Elizalde \cite{DeuEliz} gave a new interpretation of $D_{n+2}$ as the sum of the values of
the largest fixed points of all non-derangements of length $n+1$. Namely,
\begin{eqnarray*}
\sum_{k=0}^{n}(k+1)P_{n,k}=D_{n+2},
\end{eqnarray*}
which is the special case of (\ref{eqn 2.3C}) when $\mathbf{t}=(0!, 1!, 2!, \dots)$ and $m=0$.

From the previous section, one can obtain many interesting properties of $P_{n,k}$ which is left to interested readers. Furthermore, one can also
explore some new relations between $P_{n,k}$ and other classical sequences such as Bell numbers or Fibonacci numbers.

\begin{example}
By Lemma 2.1, one can derive the bivariate exponential generating function for $P_{n+k,k}$, i.e.,
\begin{eqnarray*}
P(x, y)=\sum_{n,k\geq
0}P_{n+k,k}\frac{x^n}{n!}\frac{y^k}{k!}=\frac{e^{-x}}{1-x-y}.
\end{eqnarray*}
Attracting the coefficient of $\frac{x^n}{n!}$ in $P(x, x^2)$, we have
\begin{eqnarray*}
\sum_{k=0}^{[n/2]}\binom{n}{2k}\binom{2k}{k}k!P_{n-k,k}=\sum_{k=0}^{n}(-1)^{n-k}\binom{n}{k}k!F_k,
\end{eqnarray*}
where $F_k$ is the $k$-th Fibonacci number defined by $\frac{1}{1-x-x^2}=\sum_{k\geq 0}F_kx^k$.
\end{example}

\begin{example}
In our case when $\mathbf{t}=(0!, 1!, 2!, \dots)$, (\ref{eqn 2.5}) and (\ref{eqn 2.6}) reduce to
\begin{eqnarray}
\sum_{k=0}^{n}\binom{n}{k}P_{n+m,m+k}y^{k}  &=& \sum_{k=0}^{n}(-1)^{n-k}\binom{n}{k}(m+k)!(y+1)^{k}, \label{eqn 3.1} \\
\sum_{k=0}^{n}\binom{n}{k}P_{m+k,m}y^{n-k}  &=& \sum_{k=0}^{n}\binom{n}{k}(m+k)!(y-1)^{n-k}.  \label{eqn 3.2}
\end{eqnarray}
\end{example}

It should be noted that (\ref{eqn 3.1}) and (\ref{eqn 3.2}) have close relations to the (re-normalized) Charlier
polynomials $C_n(u, v)$ \cite{Gesselb} defined by
\begin{eqnarray*}
C_n(u, v)=\sum_{k=0}^{n}\binom{n}{k}(u)_{k}v^{n-k},
\end{eqnarray*}
where $(u)_k=u(u+1)\cdots (u+k-1)$. In fact (\ref{eqn 3.1}) is equal to $\frac{(y+1)^n}{m!}C_n(m+1, -\frac{1}{y+1})$
and (\ref{eqn 3.2}) is equal to $\frac{1}{m!}C_n(m+1, y-1)$.

Recall that by (\ref{eqn 2.2C}) $P_{n,k}$ can be represented umbrally as
\begin{eqnarray*}
P_{n,k}=\mathbf{P}^{k}(\mathbf{P}-1)^{n-k},
\end{eqnarray*}
where $\mathbf{P}=\mathbf{Y_t}$ with $\mathbf{t}=(0!, 1!, 2!, \dots)$. In particular, $D_n=(\mathbf{P}-1)^{n}$ and $n!=\mathbf{P}^{n}$.
Hence, the case $y=\mathbf{P}-1$ in (\ref{eqn 3.1}) and the case $y=\mathbf{P}$ in (\ref{eqn 3.2}) generate
\begin{eqnarray*}
\sum_{k=0}^{n}\binom{n}{k}P_{n+m,m+k}D_{k}  &=& \sum_{k=0}^{n}(-1)^{n-k}\binom{n}{k}(m+k)!k!,  \\
\sum_{k=0}^{n}\binom{n}{k}P_{m+k,m}(n-k)!   &=& \sum_{k=0}^{n}\binom{n}{k}(m+k)!D_{n-k}.
\end{eqnarray*}
With the Bell umbra $\mathbf{B}$ \cite{Gesselb, Roman, RomRota}, given by $\mathbf{B}=\mathbf{Y_t}$
with $\mathbf{t}=(1, 1, 1, \dots)$. Clearly, the Bell number $B_n=\mathbf{B}^n$ and $\mathbf{B}^{n+1}=(\mathbf{B}+1)^{n}$. Then
the case $y=\mathbf{B}$ in (\ref{eqn 3.1}) and the case $y=\mathbf{B}+1$ in (\ref{eqn 3.2}) generate
\begin{eqnarray*}
\sum_{k=0}^{n}\binom{n}{k}P_{n+m,m+k}B_{k}  &=& \sum_{k=0}^{n}(-1)^{n-k}\binom{n}{k}(m+k)!B_{k+1},  \\
\sum_{k=0}^{n}\binom{n}{k}P_{m+k,m}B_{n-k+1}   &=& \sum_{k=0}^{n}\binom{n}{k}(m+k)!B_{n-k}.
\end{eqnarray*}
Using the Riordan identity \cite[P173]{Comtet},
\begin{eqnarray*}
\sum_{k=0}^{n}\binom{n}{k}(k+1)!(n+1)^{n-k}  &=& (n+1)^{n+1},
\end{eqnarray*}
the case in (\ref{eqn 3.1}) with $m=1$ and $y=-\frac{n+2}{n+1}$ and the case in (\ref{eqn 3.2}) with $m=1$ and $y=n+2$ generate respectively
\begin{eqnarray}
\sum_{k=0}^{n}(-1)^{n-k}\binom{n}{k}P_{n+1,k+1}(n+2)^{k}(n+1)^{n-k}  &=& (n+1)^{n+1},  \nonumber  \\
\sum_{k=0}^{n}\binom{n}{k}(D_k+D_{k+1})(n+2)^{n-k}  &=& (n+1)^{n+1}, \label{eqn 3.3}
\end{eqnarray}
where we use the relation $P_{k+1,1}=D_k+D_{k+1}$. By the well-known recurrence $D_{k+2}=(k+1)(D_k+D_{k+1})$ for derangement numbers $D_k$, together with $D_1=0$, (\ref{eqn 3.3}) is equivalent to
\begin{eqnarray}\label{eqn 3.4}
\sum_{k=0}^{n}\binom{n}{k}D_{k+1}(n+1)^{n-k}  &=& n^{n+1}.
\end{eqnarray}

To our best knowledge, (\ref{eqn 3.3}) and (\ref{eqn 3.4}) are the new and suprising identities analogous to the Riordan identity above.
In a forthcoming paper, using the functional digraph theory, we will give a combinatorial interpretation for a more general identity involving
the Riordan identity and (\ref{eqn 3.4}) as special cases.

\section{The special case for involutions}

In this section, we consider the special case in detail when $t_1=t_2=1$ and $t_j=0$ for $j\geq 3$. That is to study partitions of $[n+1]$ with no blocks
of sizes greater than $2$, such partitions are equivalent to involutions of $[n+1]$. Let $Q_{n,k}=A_{n,k}(\mathbf{t})$ with
$\mathbf{t}=(1, 1, 0, \dots)$, namely, $Q_{n,k}$ is the number of involutions of $[n+1]$ with the largest fixed point $k+1$.
See Table $2$ for some small values of $Q_{n,k}$. Clearly, $Q_{n,n}=I_n=\mathcal{Y}_n(1, 1, 0, \dots)$ and $Q_{n,0}=M_n=\mathcal{Y}_n(0, 1, 0, \dots)$,
where $I_n$ is the number of involutions of $[n]$, and $M_n$ is the number of involutions of $[n]$ without fixed points. It is well known
that $I_n$ and $M_n$ have the explicit formulas
\begin{eqnarray*}
I_n=\sum_{k=0}^{[n/2]}\binom{n}{2k}(2k-1)!! \hskip.5cm {\mbox{and}} \hskip.5cm
M_n=\left\{
\begin{array}{ll}
(2k-1)!!=\frac{(2k)!}{2^kk!}, & \mbox{if}\ n=2k, \\
0      , & \mbox{otherwise}.
\end{array}\right.
\end{eqnarray*}

\begin{center}
\begin{eqnarray*}
\begin{array}{c|ccccccccc}\hline
n/k & 0   & 1   & 2    & 3    & 4    &  5   & 6     & 7    & 8     \\\hline
  0 & 1   &     &      &      &      &      &       &      &         \\
  1 & 0   & 1   &      &      &      &      &       &      &        \\
  2 & 1   & 1   & 2    &      &      &      &       &      &         \\
  3 & 0   & 1   & 2    & 4    &      &      &       &      &          \\
  4 & 3   & 3   & 4    & 6    &  10  &      &       &      &         \\
  5 & 0   & 3   & 6    & 10   &  16  &  26  &       &      &         \\
  6 & 15  & 15  & 18   & 24   &  34  &  50  &  76   &      &            \\
  7 &  0  & 15  & 30   & 48   &  72  &  106 &  156  & 232  &           \\
  8 & 105 & 105 & 120  & 150  & 198  &  270 &  376  & 532  & 764       \\\hline
\end{array}
\end{eqnarray*}
Table 2. The values of $Q_{n,k}$ for $n$ and $k$ up to $8$.
\end{center}

Setting $\mathbf{t}=(1, 1, 0, \dots)$ in Lemma 2.1, one has the bivariate exponential generating function for $Q_{n+k,k}$.
\begin{eqnarray*}
Q(x, y)=\sum_{n,k\geq 0}Q_{n+k,k}\frac{x^n}{n!}\frac{y^k}{k!}=e^{y+\frac{1}{2}(x+y)^2}.
\end{eqnarray*}
Then
\begin{eqnarray*}
\sum_{n\geq 0}M_{n}\frac{x^n}{n!}=e^{\frac{1}{2}x^2} \hskip.5cm {\mbox{and}} \hskip.5cm  \sum_{k\geq 0}I_{k}\frac{y^k}{k!}=e^{y+\frac{1}{2}y^2}.
\end{eqnarray*}

Define the umbra $\mathbf{I}=\mathbf{Y_t}$ with $\mathbf{t}=(1, 1, 0, \dots)$ and $\mathbf{M}=\mathbf{Y_t}$ with $\mathbf{t}=(0, 1, 0, \dots)$, then
$\mathbf{I}=\mathbf{M}+1$, and $Q_{n,k}$ can be represented umbrally as
\begin{eqnarray}
Q_{n,k} &=& \mathbf{I}^{k}(\mathbf{I}-1)^{n-k}, \label{eqn 4.1.1}\\
Q_{n,k} &=& (\mathbf{M}+1)^{k}\mathbf{M}^{n-k}. \label{eqn 4.1.2}
\end{eqnarray}

\begin{theorem} For any integers $n,m, k\geq 0$, there hold
\begin{eqnarray}
Q_{n+k,k}&=& \sum_{j=0}^{n}(-1)^{n-j}\binom{n}{j}I_{k+j},                    \label{eqn 4.2.1} \\
Q_{n+k,k}&=& \sum_{j=0}^{k}\binom{k}{j}M_{n+j},                               \label{eqn 4.2.2}  \\
Q_{n+m+k,m+k}&=& \sum_{j=0}^{m}\binom{m}{j}Q_{n+k+j,k},                       \label{eqn 4.2.3}  \\
Q_{n+k,k}&=& \sum_{j=0}^{min\{n,k\}}\binom{k}{j}\binom{n}{j}j!I_{k-j}M_{n-j}, \label{eqn 4.2.4} \\
Q_{n+k,k}&=& \sum_{j=[\frac{n+1}{2}]}^{[\frac{n+k+1}{2}]}\frac{k!}{(n+k-2j)!}I_{n+k-2j}B(n, j),  \label{eqn 4.2.5}
\end{eqnarray}
where $B(n,j)=\frac{n!}{2^{n-j}(n-j)!(2j-n)!}$ is the Bessel number counting all the partitions of $[n]$ into $j$ blocks with the restriction of block sizes $\leq 2$.
\end{theorem}
\pf By the binomial identity, (\ref{eqn 4.2.1})-(\ref{eqn 4.2.3}) can be easily obtained by using (\ref{eqn 4.1.1}) and (\ref{eqn 4.1.2}). Attracting
the coefficient of $\frac{x^ny^k}{n!k!}$ from $Q(x,y)=e^{xy}e^{y+\frac{1}{2}y^2}e^{\frac{1}{2}x^2}$ produces (\ref{eqn 4.2.4}), and
(\ref{eqn 4.2.5}) can be derived from (\ref{eqn 4.2.4}) by shifting the index $j:=2j-n$. \qed\vskip.2cm


\begin{theorem} For any integer $n\geq 0$ and any indeterminant $y$, there hold
\begin{eqnarray}\label{eqn 4.3}
\sum_{k=0}^{n}(-1)^{n-k}\binom{n}{k}Q_{n,k}(y+1)^{k}=\sum_{k=0}^{n}\binom{n}{k}y^{k}I_{k},
\end{eqnarray}
or equivalently
\begin{eqnarray}
\sum_{k=0}^{n}\binom{n}{k}Q_{n,k}y^{k} &=& \sum_{k=0}^{n}(-1)^{n-k}\binom{n}{k}(y+1)^{k}I_{k}, \label{eqn 4.4}\\
\sum_{k=0}^{n}\binom{n}{k}Q_{n,k}y^{k} &=& \sum_{k=0}^{[n/2]}\binom{n}{2k}(2k-1)!!y^{n-2k}(y+1)^{2k}. \label{eqn 4.5}
\end{eqnarray}
\end{theorem}
\pf By (\ref{eqn 4.1.1}), we have
\begin{eqnarray*}
\sum_{k=0}^{n}(-1)^{n-k}\binom{n}{k}Q_{n,k}(y+1)^{k}
&=& \sum_{k=0}^{n}(-1)^{n-k}\binom{n}{k}(y+1)^{k}\mathbf{I}^{k}(\mathbf{I}-1)^{n-k}       \\
&=& (y\mathbf{I}+1)^n=\sum_{k=0}^{n}\binom{n}{k}y^{k}\mathbf{I}^{k}  \\
&=& \sum_{k=0}^{n}\binom{n}{k}y^{k}I_{k} , \\
\end{eqnarray*}
which proves (\ref{eqn 4.3}). Similarly, one can prove (\ref{eqn 4.4}), which can also be obtained by setting $y:=-y-1$ in (\ref{eqn 4.3}).
For (\ref{eqn 4.5}), by (\ref{eqn 4.1.2}), we have
\begin{eqnarray*}
\sum_{k=0}^{n}\binom{n}{k}Q_{n,k}y^{k}
&=& \sum_{k=0}^{n}\binom{n}{k}y^{k}(\mathbf{M}+1)^{k}\mathbf{M}^{n-k}       \\
&=& (y+(y+1)\mathbf{M})^n=\sum_{k=0}^{n}\binom{n}{k}\mathbf{M}^{k}y^{n-k}(y+1)^{k}          \\
&=& \sum_{k=0}^{n}\binom{n}{k}M_{k}y^{n-k}(y+1)^{k} , \\
&=& \sum_{k=0}^{[n/2]}\binom{n}{2k}(2k-1)!!y^{n-2k}(y+1)^{2k}.
\end{eqnarray*}

\begin{theorem} For any integers $n,m\geq 0$ and any indeterminant $y$, there hold
\begin{eqnarray}
\sum_{k=0}^{n}\binom{n}{k}Q_{m+k,m}(y+1)^{n-k}  &=& \sum_{k=0}^{n}\binom{n}{k}I_{m+k}y^{n-k}, \label{eqn 4.7}\\
\sum_{k=0}^{n}\binom{n}{k}Q_{m+k,m}B_{n-k+1}(y)  &=& y\sum_{k=0}^{n}\binom{n}{k}I_{m+k}B_{n-k}(y), \label{eqn 4.8}\\
\sum_{k=0}^{n}\binom{n}{k}\binom{y+n-k}{n-k}Q_{m+k,m}  &=& \sum_{k=0}^{n}\binom{n}{k}\binom{y}{n-k}I_{m+k}. \label{eqn 4.9}
\end{eqnarray}
\end{theorem}
\pf The special case in (\ref{eqn 2.6}) with $\mathbf{t}=(1, 1, 0, \dots)$ and $y:=y+1$ generates (\ref{eqn 4.7}).
For (\ref{eqn 4.8}), define a linear (invertible) transformation
\begin{eqnarray*}
L_1(y^k)=B_k(y),\ \ \ (k=0,1,2,\dots),
\end{eqnarray*}
where $B_k(y)=\mathcal{Y}_{k}\big(y, y, y,\dots\big)$ is the Bell polynomial satisfying the relation
\begin{eqnarray*}
B_{n+1}(y)=y\sum_{k=0}^{n}\binom{n}{k}B_k(y).
\end{eqnarray*}
Then we have
\begin{eqnarray*}
yL_1((y+1)^n)=y\sum_{k=0}^{n}\binom{n}{k}L_1(y^k)=y\sum_{k=0}^{n}\binom{n}{k}B_k(y)=B_{n+1}(y).
\end{eqnarray*}
Hence (\ref{eqn 4.8}) follows by acting $yL_1$ on the two sides of (\ref{eqn 4.7}).

For (\ref{eqn 4.9}), similarly, define another linear transformation
\begin{eqnarray*}
L_2(y^k)=\binom{y}{k},\ \ \ (k=0,1,2,\dots).
\end{eqnarray*}
By the Vandermonde's convolution identity, we have
\begin{eqnarray*}
L_2((y+1)^n)=\sum_{k=0}^{n}\binom{n}{k}L_2(y^k)=\binom{y+n}{n}.
\end{eqnarray*}
Then acting $L_2$ on the two sides of (\ref{eqn 4.7}) leads to (\ref{eqn 4.9}). \qed\vskip0.2cm

\begin{corollary} For any integers $n,m \geq 0$, there hold
\begin{eqnarray}
\sum_{k=0}^{n}\binom{n}{k} Q_{m+k,m}(n-k)!         &=& \sum_{k=0}^{n}\binom{n}{k}I_{m+k}D_{n-k},    \label{eqn 4.10} \\
\sum_{k=0}^{n}\binom{n}{k}Q_{m+k,m}I_{n-k}         &=& \sum_{k=0}^{n}\binom{n}{k}I_{m+k}M_{n-k},     \label{eqn 4.11} \\
\sum_{k=0}^{n}\binom{n}{k} Q_{m+k,m}B_{n-k+1}      &=& \sum_{k=0}^{n}\binom{n}{k}I_{m+k}B_{n-k},      \label{eqn 4.12}
\end{eqnarray}
\end{corollary}
\pf Setting $y=\mathbf{D}$, $y=\mathbf{M}$ and $y=\mathbf{B}$ in (\ref{eqn 4.7}) produces (\ref{eqn 4.10})-(\ref{eqn 4.12}) respectively.
\qed\vskip0.2cm

The special cases in (\ref{eqn 2.3C}) and (\ref{eqn 2.3D}) with $\mathbf{t}=(1, 1, 0, \dots)$ and $m=0$ generate

\begin{theorem} For any integer $n\geq 0$, there hold
\begin{eqnarray*}
\sum_{k=0}^{n}(k+1)Q_{n,k}  &=& M_{n+2}-I_{n+2}+(n+2)I_{n+1},   \\
\sum_{k=0}^{n}(n-k+1)Q_{n,k}&=& I_{n+2}-(n+2)M_{n+1}-M_{n+2} .
\end{eqnarray*}
\end{theorem}

\section{The special case for labeled forests}

In this section, we consider the special case when $t_j=j^{j-1}$ for $j\geq 1$. That is to assign a (rooted and labeled)
tree structure to each block of partitions of $[n+1]$, such partitions with weight $\mathbf{t}=(1^0, 2^1, 3^2, \dots)$ are equivalent to
labeled forests on $[n+1]$. Let $L_{n,k}=A_{n,k}(\mathbf{t})$ with
$\mathbf{t}=(1^0, 2^1, 3^2, \dots)$, namely, $L_{n,k}$ is the number of labeled forests on $[n+1]$ with the largest singleton tree labeled by $k+1$.
A {\em singleton tree} is a labeled tree with exactly one point. Clearly,
$L_{n,n}=\mathcal{Y}_n(1^0, 2^1, 3^2, \dots)=(n+1)^{n-1}$ and $L_{n,0}=\mathcal{Y}_n(0, 2^1, 3^2, \dots)$, where $L_{n,0}$
is also the number of labeled forests on $[n]$ with no singleton trees. See Table $3$ for some small values of $L_{n,k}$.

\begin{center}
\begin{eqnarray*}
\begin{array}{c|ccccccccc}\hline
n/k & 0     & 1     & 2     & 3     & 4      &  5    & 6          \\\hline
  0 & 1     &       &       &       &        &       &            \\
  1 & 0     & 1     &       &       &        &       &           \\
  2 & 2     & 2     & 3     &       &        &       &            \\
  3 & 9     & 11    & 13    & 16    &        &       &             \\
  4 & 76    & 85    & 96    & 109   & 125    &       &            \\
  5 & 805   & 881   & 966   & 1062  & 1171   & 1296  &             \\
  6 & 10626 & 11431 & 12312 & 13278 & 14340  & 15511 & 16807          \\\hline
\end{array}
\end{eqnarray*}
Table 3. The values of $L_{n,k}$ for $n$ and $k$ up to $6$.
\end{center}

Setting $\mathbf{t}=(1^0, 2^1, 3^2, \dots)$ in Lemma 2.1, and using the identity \cite[P174]{Comtet}
\begin{eqnarray*}
\exp{\big(\sum_{j\geq 1}j^{j-1}\frac{z^j}{j!}\big)}
=\sum_{j\geq 0}(j+1)^{j-1}\frac{z^j}{j!} ,
\end{eqnarray*}
one has the bivariate exponential generating function for $L_{n+k,k}$,
\begin{eqnarray*}
L(x, y)=\sum_{n,k\geq 0}L_{n+k,k}\frac{x^n}{n!}\frac{y^k}{k!}
=e^{-x}\big(\sum_{j\geq 0}(j+1)^{j-1}\frac{(x+y)^j}{j!}\big).
\end{eqnarray*}

Define the umbra $\mathbf{L}=\mathbf{Y_t}$ with $\mathbf{t}=(1^0, 2^1, 3^2, \dots)$, then $L_{n,k}$ can be represented umbrally as
\begin{eqnarray}
L_{n,k} &=& \mathbf{L}^{k}(\mathbf{L}-1)^{n-k}, \label{eqn 5.1.1}
\end{eqnarray}

Similar to the Section 4, using (\ref{eqn 5.1.1}) one can derive the corresponding results for $L_{n+k,k}$, the details are left to readers.

\begin{theorem} For any integers $n,m, k\geq 0$, there hold
\begin{eqnarray*}
L_{n+k,k} &=& \sum_{j=0}^{n}(-1)^{n-j}\binom{n}{j}(k+j+1)^{k+j-1},                    \label{eqn 5.2.1} \\
L_{n+m+k,m+k}&=& \sum_{j=0}^{m}\binom{m}{j}L_{n+k+j,k},                       \label{eqn 5.2.3} \\
\sum_{j=0}^{n}\binom{j+\lambda-1}{j}L_{n+m,m+j} &=& \sum_{j=0}^{n}(-1)^{n-j}\binom{n+\lambda}{j}(m+j+1)^{m+j-1}.
\end{eqnarray*}
\end{theorem}

\begin{theorem} For any integer $n\geq 0$ and any indeterminant $y$, there hold
\begin{eqnarray*}\label{eqn 5.3}
\sum_{k=0}^{n}(-1)^{n-k}\binom{n}{k}L_{n,k}(y+1)^{k}=\sum_{k=0}^{n}\binom{n}{k}(k+1)^{k-1}y^{k},
\end{eqnarray*}
or equivalently
\begin{eqnarray*}
\sum_{k=0}^{n}\binom{n}{k}L_{n,k}y^{k} &=& \sum_{k=0}^{n}(-1)^{n-k}\binom{n}{k}(k+1)^{k-1}(y+1)^{k}. \label{eqn 5.4}
\end{eqnarray*}
\end{theorem}

\begin{theorem} For any integers $n,m\geq 0$ and any indeterminant $y$, there hold
\begin{eqnarray*}
\sum_{k=0}^{n}\binom{n}{k}L_{m+k,m}(y+1)^{n-k}  &=& \sum_{k=0}^{n}\binom{n}{k}(m+k+1)^{m+k-1}y^{n-k}, \label{eqn 5.7}\\
\sum_{k=0}^{n}\binom{n}{k}L_{m+k,m}B_{n-k+1}(y)  &=& y\sum_{k=0}^{n}\binom{n}{k}(m+k+1)^{m+k-1}B_{n-k}(y), \label{eqn 5.8}\\
\sum_{k=0}^{n}\binom{n}{k}\binom{y+n-k}{n-k}L_{m+k,m}  &=& \sum_{k=0}^{n}\binom{n}{k}\binom{y}{n-k}(m+k+1)^{m+k-1}. \label{eqn 5.9}
\end{eqnarray*}
\end{theorem}

\begin{corollary} For any integers $n,m \geq 0$, there hold
\begin{eqnarray*}
\sum_{k=0}^{n}\binom{n}{k} L_{m+k,m}(n-k)!         &=& \sum_{k=0}^{n}\binom{n}{k}(m+k+1)^{m+k-1}D_{n-k},    \label{eqn 5.10} \\
\sum_{k=0}^{n}\binom{n}{k} L_{m+k,m}B_{n-k+1}      &=& \sum_{k=0}^{n}\binom{n}{k}(m+k+1)^{m+k-1}B_{n-k},      \label{eqn 5.12}
\end{eqnarray*}
\end{corollary}

\section*{Acknowledgements} The authors are grateful to the
anonymous referees for the helpful suggestions and comments. The
work was supported by The National Science Foundation of China
(Grant No. 10801020 and 70971014) and supported by the Fundamental
Research Funds for the Central Universities (Grant No. 2009GN070).


\end{document}